\documentclass[12pt]{article}

\usepackage{centernot}
\usepackage{bbm}
\usepackage{booktabs}
\usepackage{float}
\usepackage[toc,page]{appendix}
\usepackage{pdfpages}
\newcounter{note}

\usepackage[expansion=false]{microtype}
\usepackage[reqno]{amsmath}
\usepackage{fixmath}
\usepackage{amssymb}
\usepackage{mathtools}
\usepackage{graphicx}
\usepackage{amsthm} 
\usepackage{enumerate}  
\usepackage{url} 
\usepackage{makeidx} 
\usepackage{tikz}
\usetikzlibrary{decorations.markings}
\usepackage{array}
\usepackage{float}
\usepackage{wrapfig}
\usepackage{graphicx}
\usepackage{ifthen, amsfonts, color, pgf}
\usepackage{hyperref}

\newtheoremstyle{plainsl}%
{\topsep}
{\topsep}
{\slshape} 
{}
{\normalfont\bfseries}
{.}
{ }
{}

\swapnumbers

{\theoremstyle{plainsl}
	\newtheorem{theorem}{Theorem}[section]
	
	}
\newtheorem*{theorem*}{Theorem}

\newtheorem{lem}[theorem]{Lemma}
\newtheorem*{lem*}{Lemma}

\newtheorem{cor}[theorem]{Corollary}

\newtheorem*{defn*}{Definition}

{\theoremstyle{remark}
	
	}

\newcommand\cref[1]{Corollary~\ref{cor:#1}}

\renewcommand\proof{\noindent\textsl{Proof. }}

\newcommand\sqr[2]{{\vbox{\hrule height.#2pt
			\hbox{\vrule width.#2pt height#1pt \kern#1pt
				\vrule width.#2pt}\hrule height.#2pt}}}
\renewcommand\qed{%
	\ifmmode\eqno\sqr53
	\else\nolinebreak\ \hfill\sqr53\medbreak\fi}

\numberwithin{equation}{section}

\newcommand\cx{{\mathbb C}}

\newcommand\onesum{{\uplus}}

\usepackage{caption}
\usepackage{multirow}
\usepackage{amsmath}

\tikzset{every node/.style={shape=circle, fill=black, inner sep=0pt, minimum size=5pt}}

\makeindex

\title{Limbs and Cospectral Vertices in Trees}
\author{Xiaojing Wang and Karen Yeats}

\begin{document}

\maketitle

\begin{abstract}
	We generalize Schwenk's result that almost all trees contain any given limb to trees with positive integer vertex weights. The concept of characteristic polynomial is extended to such weighted trees and we prove that the proportion of $n$-vertex weighted trees that is weighted cospectral to another $n$-vertex weighted tree approaches $1$ as $n$ approaches infinity. We also prove that, for any integer $k\ge 2$, the proportion of $n$-vertex trees containing $k$ non-similar cospectral vertices approaches $1$ as $n$ approaches infinity.
\end{abstract}

\allowdisplaybreaks

\section{Introduction}

Let $G$ be a graph and suppose $A$ is the adjacency matrix of $G$. The \emph{characteristic polynomial} $\phi(G, x)$ of $G$ is defined as $$\phi(G, x)=\det(xI-A).$$ When there is no ambiguity, we use $\phi(G)$ to denote the characteristic polynomial of $G$. Two graphs are \emph{cospectral} if they have the same characteristic polynomial. Whether a graph is cospectral to another graph is considered in a variety of capacities in algebraic graph theory. van Dam and Haemer's \cite{vDH2003} survey pointed out the point graph and the line graph of a partial linear space are cospectral if the number of points on each line equals the number of lines a point is on. Godsil and MCKay\cite{GM1982} found a local switching construction for cospectral graphs, but the local switching can only be done for graphs with a partition satisfying a set of specific requirements. Out of the known constructions of cospectral graphs, one notable result was due to Schwenk \cite{Schwenk1973}, where he proved that as the integer $n$ approaches infinity, the proportion of $n$-vertex trees that is cospectral to another tree approaches $1$.

To prove this result, Schwenk used a particular type of subtree, which he called a limb \cite{Schwenk1973}. For any tree $T$ and some vertex $v\in V(T)$, a \emph{branch} at $v$ is a maximal subtree of $T$ where $v$ has degree $1$. A \emph{limb} at $v$ is the rooted subtree formed by a collection of branches at $v$, and $v$ is the root of the limb. Note the tree $T$ is not rooted, but the limb is.

Schwenk proved that, given an $\ell$-vertex rooted tree $L$, the number of $n$-vertex trees that contains $L$ as a limb is a constant that only depends on $n$ and $\ell$ \cite[Theorem 4]{Schwenk1973}. The following theorem is a rewording of his result.

\begin{theorem}[\cite{Schwenk1973}]
	If $R$ and $S$ are two rooted trees with $\ell$ vertices, and $r_n$ (or $s_n$) is the number of $n$-vertex trees which do not have $R$ (or $S$) as a limb, then $r_n=s_n$.
	\label{lreplacedsamect}
\end{theorem}

This theorem is interesting because it shows the structure of the limb is irrelevant when counting the number of trees over $n$ vertices with a specific limb. Using the structure-less property in Theorem \ref{lreplacedsamect}, any $\ell$ vertex limb can be treated at $K_{1, \ell-1}$ for the purpose of counting the number of a specified $\ell$-vertex limb in $n$-vertex trees. Schwenk used this fact to derive the following recursive equation that the generating function $S(x)$ of trees without a specified $\ell$-vertex limb must satisfy \cite[Theorem 5]{Schwenk1973}.

$$S(x)=(x-x^\ell)\sum_{i\ge1}\frac{1}{i}S(x^i)$$

It is worth noting that McAvaney \cite{McAvaney1974} and Lu \cite{Lu1996} independently derived this recursive equation using more direct enumerative arguments. 

By performing an asymptotic analysis, Schwenk showed the radius of convergence of $S(x)$ is greater than the radius of convergence of the generating function of trees \cite{Schwenk1973}. This implies that the proportion of $n$-vertex trees containing a specified limb approaches $1$ as $n$ approaches infinity. Informally, we say that almost all trees has a cospectral mate. In this paper, we use this meaning of almost all in our discussions.

Schwenk observed there exists a pair of limbs on $9$ vertices, such that if one is replaced by the other in a tree, the characteristic polynomial of the tree remains unchanged \cite{Schwenk1973}. As the integer $n$ approaches infinity, the proportion of $n$-vertex trees containing at least one of the pair of limbs approaches $1$. the following theorem, which is a rewording of Schwenk's original result, was proven.

\begin{theorem}[\cite{Schwenk1973}]
	For any $0\le p\le 1$, there exists a positive integer $n$, such that for any positive integer $N\ge n$, the proportion of $N$-vertex trees that is cospectral to another tree is greater than $p$.
	\label{treecosmate}
\end{theorem}

\subsection{Limbs in Rooted Trees}

In order to derive the generating function of tress without a specified $l$-vertex limb, McAvaney discussed how to find the generating function of rooted trees without a specified $l$-vertex limb \cite{McAvaney1974}. The generating function he found implies the following extension to rooted trees of Schwenk's result.

\begin{theorem}[\cite{McAvaney1974}]
	Let $\mathcal{S}$ be the set of rooted trees without some $\ell$-vertex rooted tree $L$ as a limb, and let $S(x)=\sum_{i=0}^\infty S_ix^i$ be its generating function, then \[S(x)=(x-x^\ell)\left(\sum_{i=1}^\infty\frac{1}{i}S(x^i)\right).\]
	\label{McAvaneygen}
\end{theorem}

The variable $\ell$ is the only quantity about $L$ in this generating function $S(x)$ in Theorem \ref{McAvaneygen}, so the number of rooted trees over $n$ vertices without $L$ as a limb is a fixed constant that only depends on $n$ and $\ell$, but not the structure of $L$. In other words, this derivation process implies the following result, which is a rooted version of Theorem \ref{lreplacedsamect}.

\begin{theorem}
	Let $L_1$ and $L_2$ be two distinct rooted trees with $\ell$ vertices, where $\ell\ge3$. Let $\ell_{1, n}$ and $\ell_{2, n}$ denote the number of $n$-vertex rooted trees without $L_1$ or $L_2$ as a limb, respectively. Then $\ell_{1, n}=\ell_{2, n}$.
	\label{rtsamenum}
\end{theorem}

McAvaney showed how Theorem \ref{McAvaneygen} coinsides with Schwenk's result in \cite{Schwenk1973}, so that Schwenk's result can be proven both combinatorically and enumeratively. Moreover, McAvaney's result is interesting because now we know the number of $n$-vertex rooted tress without $L$ as a limb does not depend on the structure of $L$.

\subsection{Weighted Trees}

Weighted trees, both rooted and unrooted, are of our interest in this paper. We define \emph{weighted rooted trees} as rooted trees where each vertex is assigned a positive integer weight. Similarly, a \emph{weighted unrooted tree}, or simply \emph{weighted tree} when there is no ambiguity, is defined as a tree where each vertex is assigned a positive integer weight, but with no distinguished root vertex.

We are interested in considering Schwenk's results in a weighted trees context. However, observe that the concept of characteristic polynomial and cospectrality were undefined for weighted trees. Therefore, to further extend Schwenk's results, we need to extend the definition of characteristic polynomial to weighted trees. For a weighted tree $W$, we use $A(W)$ to denote the adjacency matrix of $W$. Suppose $I^*(x)$ is an $n\times n$ diagonal matrix, such that its rows and columns are each indexed by the vertices of $W$, and the $ii$-entry of $I^*(x)$ is $x^{w(i)}$, where $w(i)$ is the weight of the vertex $i$ and $x$ an indeterminate. Let $A(W)$ be the adjacency matrix of $W$. The \emph{weighted characteristic polynomial}\index{weighted characteristic polynomial} $\phi^*(W, x)$ of $W$ is defined as \[\phi^*(W, x)=\det(M(x)-A(W)).\] When there is no ambiguity, we drop the indeterminate variable and use $\phi^*(W)=\det(I^*-A(W))$ to denote the weighted characteristic polynomial of $W$. If two weighted graphs have the same weighted characteristic polynomial, we say that they are \emph{weighted cospectral}.

In this paper, we extend several of Schwenk's results and McAvaney's results, such as Theorem \ref{lreplacedsamect} and Theorem \ref{treecosmate}, in a variety of contexts. We show that the number of weighted rooted trees over $n$ vertices containing a specified limb is also independent from the structure of the limb. We deduce that, for a given rooted tree $L$, the proportion of $n$-vertex weighted rooted trees and $n$-vertex weighted trees with $L$ as a limb approaches $1$ as $n$ approaches infinity. We also consider a special type of limb that takes all branches at a given vertex, and prove similar results apply to this type of limbs as well. A property analogous to cospectrality is discussed for weighted rooted trees. Furthermore, using Schwenk's result in Theorem \ref{treecosmate}, we give a construction for trees with $k$ cospectral vertices, for any integer $k\ge2$. 


\section{Limbs in Weighted Rooted Trees}\label{sec weighted limbs}

In this section, we extend Schwenk's limb replacement result and asymptotic analysis to weighted rooted trees.

We use the generating function argument that Schwenk applied to perform his analysis. To utilize this argument, we need to first derive a recursion for the generating function of the set of weighted rooted trees. Let $\mathcal{T_W}$ be the set of all weighted rooted trees, and let $T_W(x)$ be its generating function, where a weight of a tree is the sum of the weights of all its vertices. We start by deriving a recursive formula for $T_W(x)$. Like for the rooted trees, given a weighted rooted tree $T_W$, we use $F(T_W)$ to denote the set of weighted rooted trees obtained by first deleting the root of $T_W$ and then letting the vertex that is the neighbor of the original root serve as the new root of each component, and where the weight on each remaining vertex is unchanged. Similarly to unweighted rooted trees, consider what happens to a weighted rooted tree $T_W$ if its root $r$ with weight $w_r$ is deleted. The $F(T_W)$ is a forest where each component of which is an element of $\mathcal{T}_R$. Use $\mathcal{R}$ to denote the set of all the possible choices for the root. Using $\mathbf{MSET}$ to denote the multiset operator, we obtain the following combinatorial isomorphism using the notation in \cite{FS2009}. (Refer to \cite{FS2009}[I.2, page 26] for more details about the $\mathbf{MSET}$ operator.)

$$\mathcal{T_W}=\mathcal{R}\times\mathbf{MSET}(\mathcal{T_W}).$$

Note the weight of $r$ can be any positive integer, therefore, we get

\begin{equation}\label{eq weighted spec}
T_W(x)=\left(\sum_{j=1}^\infty x^j\right)\exp\left(\sum_{i=1}^\infty\frac{1}{i}T(x^i)\right)=\frac{x}{1-x}\exp\left(\sum_{i=1}^\infty\frac{1}{i}T_W(x^i)\right).
\end{equation}

Observe that rooted trees can be considered as a specific case of weighted rooted trees where all vertices have weight $1$.





Analogous to Schwenk's proof, we start by considering the radius of convergence of $T_W(x)$, denoted as $\alpha_T$.

\begin{lem}
	The radius of convergence of $T_W(x)$ is at least $\frac{1}{16}$.
\end{lem}

\proof Define the sequence $\{T_n\}$ such that $$T_W(x)=\sum_{n\geq 1}T_nx^n.$$ Let $$A(x)=\frac{T_W(x)}{x}$$ and $$B(x)=\sum_{i\geq 1}\frac{1}{i}T_W(x^i).$$ Note that, by \eqref{eq weighted spec}, $A(x)=\frac{1}{1-x}\exp{B(x)}$. Moreover, we have 
\begin{eqnarray*}
	\frac{d}{dx}A(x)&=&\frac{1}{1-x}\exp\left(\sum_{i=1}^\infty\frac{1}{i}T_W(x^i)\right)\cdot\frac{d}{dx}B(x)\\&&+\frac{1}{(1-x)^2}\exp\left(\sum_{i=1}^\infty\frac{1}{i}T_W(x^i)\right)\\
	&=&A(x)\cdot\frac{d}{dx}B(x)+\frac{1}{1-x}A(x).
\end{eqnarray*}

Therefore, using the notation $\langle x^m, FGx) \rangle$ for the coefficient of $x^m$ in $G(x)$, 

\begin{eqnarray*}
	T_{n+1}&=&\frac{1}{n}\langle x^{n-1}, \frac{d}{dx}A(x)\rangle\\
	&=&\frac{1}{n}\langle x^{n-1}, A(x)\cdot\frac{d}{dx}B(x)\rangle+\frac{1}{n}\sum_{i=0}^{n-1}\langle x^{n-1-i}, A(x)\rangle\\
	&=&\frac{1}{n}\sum_{i=1}^n\left(\sum_{d|i}dT_d\right)T_{n-i+1}+\frac{1}{n}\sum_{i=1}^nT_{n-i+1}\\
	&=&\frac{1}{n}\sum_{i=1}^n\left(iT_i\cdot\sum_{i\le ci\le n}T_{n-ci+1}+T_{n-i+1}\right)\\
\end{eqnarray*}

For every weighted rooted tree of weight $n$, by adding $1$ to the weight of its root, we can create a distinct weighted rooted tree corresponding to it with weight $n+1$. Therefore, the sequence $\{T_n\}$ is weakly increasing, and so for any $c$ such that $i\le ci\le n$, we have $T_{n-ci+1}\le T_{n-i+1}$. Then $\sum_{i\le ci\le n}T_{n-ci+1}\le nT_{n-i+1}$. Additionally, since $T_1=1$, it is clear that $T_i\ge 1$ for any $i\ge 1$, so $\frac{T_{n-i+1}}{iT_i}\le T_{n-i+1}$ for any $1\le i\le n$. Thus,

\begin{eqnarray*}
	&=&\frac{1}{n}\sum_{i=1}^niT_i\cdot\left(\sum_{i\le ci\le n}T_{n-ci+1}+\frac{T_{n-i+1}}{iT_i}\right)\\
	&\le&\frac{1}{n}\sum_{i=1}^niT_i\cdot\frac{n+i}{i}T_{n-i+1}\\
	&\le&2\sum_{i=1}^nT_iT_{n-i+1}.
\end{eqnarray*}

Now, we define a power series $f(x)=\sum_{i=1}f_ix^i$, where $f_1=1$ and $f_n=2\sum_{i=1}^nf_if_{n-i-1}$. Since $f(x)$ bounds $T_W(x)$ above, the radius of convergence of $f(x)$ is a lower bound of the radius of convergence of $T_W(x)$. Moreover, $\langle x^{n+1}, (f(x))^2\rangle=\sum_{i=1}^nf_if_{n-i+1}=\frac{1}{2}f_{n+1}$ for all $n\ge2$, and $\langle x, (f(x))^2\rangle=0$. Then we have $$(f(x))^2-\frac{1}{2}f(x)+x=0.$$

By the quadratic formula and the fact that $f(0)=0$, we get $$f(x)=\frac{1}{4}\left(1-\frac{1}{\sqrt{1-16x}}\right).$$ The radius of convergence of $f(x)$ is $\frac{1}{16}$, so the radius of convergence of $T_W(x)$ is at least $\frac{1}{16}$. \qed

The following lemma shows that $T_W(x)$ is bounded when $x=\alpha_T$, where $\alpha_T$ the radius of convergence of $T_W(x)$. We will use this to compare $\alpha_T$ and the radius of convergence of the generating series of weighted trees with a forbidden limb. The general approach used to prove this lemma is described in Section 9.5 of Harary and Palmer \cite{HP1973}.

\begin{lem}
	The series $T_W(x)$ satisfies $T_W(\alpha_T)=1$.
	\label{radT1}
\end{lem}

\proof Define a multivariate function: for $x, y\in\cx$, $$G(x, y)=\frac{x}{1-x}\exp{\left(y+\sum_{i=2}^\infty\frac{1}{i}T_W(x^i)\right)}-y.$$

By the Implicit Function Theorem, $y=T_W(x)$ is the unique analytic solution of $G(x, y)=0$. Moreover, it has a singularity at $x=\alpha_T$, and $G(\alpha_T, T_W(\alpha_T))=0$.

Therefore, $$\frac{\partial G}{\partial y}(\alpha_T, T_W(\alpha_T))=G(\alpha_T, T_W(\alpha_T))+T_W(x)-1=T_W(\alpha_T)-1=0,$$ due to the singularity. Thus $T_W(\alpha_T)=1$. \qed

Now we've established the groundwork to discuss weighted rooted trees, let's define what a limb is in this case. Let $T$ be a weighted rooted tree and let $v\in V(T)$. A \emph{branch} $B$ at $v$ is a weighted rooted subtree of $T$, such that $v$ is the root of $B$, $v$ has degree $1$ in $B$, and if $v'$ is a children of $v$ and $v'$ is in $T$, then all descendants of $v'$ are in $T$.
A \emph{limb} at $v$ is a weighted rooted subtree of $T$ with root $v$ that consists of a collection of branches at $v$.

Let $\mathcal{S}_W$ be the set of all weighted rooted trees without a specific weighted limb $L$. We derive its generating series $S_W(x)$ by an argument similar to Schwenk's \cite{Schwenk1973}.
Suppose the weight of $L$ is $\ell$, and the weight of its root vertex is $w$. If $L$ is too small, $\mathcal{S}_W$ would be trivial. So, assume $L$ contains at least three vertices, which implies $\ell\ge3$. 

We construct the generating function $S_W(x)$ using the recursive relationship among the trees in the set $\mathcal{S}_W$.
Given a weighted rooted tree $T_W$ in $\mathcal{S}_W$, observe that $F(T_W)$ is a multiset of weighted rooted trees in $S_W$ that does not contain $F(L)$ as a submultiset. Assume the root of $L$ has weight $w$, then the weight of $F(L)$ is $\ell-w$. The generating function $S_W(x)$ is the difference between the generating function for weighted rooted trees with a root of any weight attached to any forest where each component is an element of $\mathcal{S}_W$ and the generating function for weighted rooted trees $U$ with root of weight $w$ while $F(U)$ is a superset of $F(L)$. If we use $\mathcal{F(L)}$ to denote the set containing exactly the forest $F(L)$, then we obtain the following recursion.

$$\mathcal{S_W}=(\mathcal{R}\times\mathbf{MSET}(\mathcal{S_W})) - (\mathcal{R}\times\mathcal{F(L)}\times\mathbf{MSET}(\mathcal{S_W})).$$

The set subtraction is allowable in this particular combinatorial specification as we can view $\mathcal{F(L)}\times\mathbf{MSET}(\mathcal{S_W})$ as a subset of $\mathbf{MSET}(\mathcal{S_W})$. The injective map giving this containment can be defined as follows: Each element of $F(L)$ is in $\mathcal{S_W}$. Since every element of the cartesian product has the same first coordinate (namely $F(L)$), the only information $\mathcal{F(L)}\times\mathbf{MSET}(\mathcal{S_W})$ carries is how many of each element of $\mathcal{S_W}$ is in the multiset. So, we obtain the injective map as follows: given an element of the cartesian product, take the union of its first and second coordinates.

Consequently, the generating function for weighted rooted trees without a specific weighted limb $L$ satisfies the following.

\begin{align*}
	S(x)&=\frac{x}{1-x}\exp\left(\sum_{i=1}^\infty\frac{1}{i}S(x^i)\right)-x^w\cdot x^{\ell-w}\exp\left(\sum_{i=1}^\infty\frac{1}{i}S(x^i)\right)\\
	&=\left(\frac{x}{1-x}-x^\ell\right)\exp\left(\sum_{i=1}^\infty\frac{1}{i}S(x^i)\right)
\end{align*}



Since the number of weighted rooted trees without a specific limb increases as a function of the number of vertices, it is easy to see that the radius of convergence of $S$, denoted as $\alpha_S$, is finite. Meanwhile, by definition, $\langle x^i, S_W(x)\rangle\le\langle x^i, T_W(x)\rangle$ for any non-negative integer $i$, so $\alpha_S\ge\alpha_T>0$. To further these two quantities, in the following lemma, we prove $S(\alpha_S)=1$.

\begin{lem}
	The series $S_W(x)$ satisfies $S_W(\alpha_S)=1$.
	\label{radS2}
\end{lem}

\proof Define a multivariate function: for $x, y\in\cx$, $$F(x, y)=\left(\frac{x}{1-x}-x^l\right)\exp{\left(y+\sum_{i=2}^\infty\frac{1}{i}S_W(x^i)\right)}-y.$$

By the Implicit Function Theorem, $y=S_W(x)$ is the unique analytic solution of $F(x, y)=0$. Moreover, it has a singularity at $x=\alpha_S$, and $F(\alpha_S, S_W(\alpha_S))=0$.

Therefore, $$\frac{\partial F}{\partial y}(\alpha_S, S_W(\alpha_S))=F(\alpha_S, S_W(\alpha_S))+S_W(x)-1=S_W(\alpha_S)-1=0,$$ due to the singularity. Thus $S_W(\alpha_S)=1$. \qed

Now we are ready to prove $\alpha_S>\alpha_T$.

\begin{theorem}
	The radius of convergence of $S_W(x)$ is greater than the radius of convergence of $T_W(x)$.
	\label{swtwrds}
\end{theorem}

\proof By definition, $S_W(x)$ is coefficient-wise bounded above by $T_W(x)$, which implies $\alpha_S\ge\alpha_T$, and $S_W(x)\le T_W(x)$ for any $x>0$. Moreover, since $\mathcal{T_W}$ contains all weighted rooted trees with weight $\ell$ while $L\not\in\mathcal{S_W}$, the coefficient of $x^l$ in $S_W$ is strictly less than the coefficient of that in $T_W$. Meanwhile, Lemma \ref{radS2} and Lemma \ref{radT1} imply $S_W(\alpha_S)=T_W(\alpha_T)$. Therefore, $\alpha_S\ne\alpha_T$, so $\alpha_S>\alpha_T$. \qed

A direct consequence of the theorem is the following result analogous to Schwenk's \cite{Schwenk1973}.

\begin{cor}
	For any given weighted rooted tree $L$, the proportion of $n$-vertex weighted rooted trees with $L$ as a limb approaches $1$ as $n$ approaches infinity.
	\label{haslimbwr}
\end{cor}

\proof By Theorem \ref{swtwrds}, $S_W(x)$ converges on a larger disk than $T_W(x)$. In other words, the coefficients of $T_W(x)$ have a larger order of growth that those of $S_W(x)$. Since $S_W(x)$ is the generated series of weighted rooted trees without $L$ as a limb, we can conclude that for any rooted weighted tree $L$, the proportion of $n$-vertex weighted rooted trees with $L$ as a limb approaches $1$ as $n$ approaches infinity. \qed

Informally, we say that almost all weighted rooted trees have any weighted rooted tree $L$ as a limb.

\section{Limbs in Weighted Unrooted Trees}
\label{sec weighted trees}

Using the underlying relationships between corresponding rooted and unrooted structures, we extend the discussion from Section \ref{sec weighted limbs} to weighted trees that are not rooted. 

Let $\mathcal{W}^\bullet$ be the set of all pairs $(T, v)$, where $T$ is a weighted tree and $v$ is a vertex in $G$, then $$\mathcal{W}^\bullet\cong\mathcal{T_W},$$ where $\mathcal{T_W}$ is the set of weighted rooted trees, as in the previous section. Similarly, let $\mathcal{W}^{\bullet-\bullet}$ be the set of all pairs $(T, e)$, where $T$ is a weighted tree and $e$ is an edge in $T$. Note we can view $T$ as a tree formed by attaching a weighted rooted tree on each end of $e$. Let $\mathbf{SET}_2(\mathcal{T_W})$ denote the set of all subsets of size $2$ of $\mathcal{T_W}$, then $$\mathcal{W}^{\bullet-\bullet}\cong\mathbf{SET}_2(\mathcal{T_W}).$$

On the other hand, define $D(\mathcal{T_W}\times\mathcal{T_W}):=\{(T, T)|T\in\mathcal{T_W}\}$, then the generating function of $D(\mathcal{T_W}\times\mathcal{T_W})$ is $T_W(x^2)$. Notice that $\mathbf{SET}_2(\mathcal{T_W})$ satisfies the combinatorial equivalence $$\mathbf{SET}_2(\mathcal{T_W})+\mathbf{SET}_2(\mathcal{T_W})=(\mathcal{T_W}\times\mathcal{T_W})+D(\mathcal{T_W}\times\mathcal{T_W}).$$

Consequently, the generating function for $\mathbf{SET}_2(\mathcal{T_W})$, or $\mathcal{W}^{\bullet-\bullet}$ is $$\frac{1}{2}((T_W(x))^2+T_W(x^2)).$$

Meanwhile, observe that the dissimilarity theorem for trees still holds for weighted trees as we simply ignore the weights for the constructions of the dissimilarity theorem (The proof for the dissimilarity theorem for trees can be found in \cite{FS2009}[VII.26, page 481].) Specifically the dissimilarity theorem says that if we let $\mathcal{W}$ be the set of all weighted trees, then\footnote{The proof of the dissimilarity theorem relies on the fact that given a tree, every maximum length path in that tree has the same center, either a vertex or an edge depending on parity, then the rootings at edges or vertices on the left correspond either to the first term on the right by rooting at the center edge or vertex, or to the second term by cutting at the rooting edge or rooting vertex in the direction of the center.} $$\mathcal{W}^\bullet+\mathcal{W}^{\bullet-\bullet}\cong\mathcal{W}+(\mathcal{W}^\bullet\times\mathcal{W}^\bullet).$$ Let $W(x)$ be the generating function of $\mathcal{W}$. Based on the combinatorial equivalence above, we get $$W(x)=T_W(x)+\frac{1}{2}(T_W(x^2)-(T_W(x))^2).$$

There is a similar relationship between weighted unrooted trees and weighted rooted trees without a specified limb.  
Let $\mathcal{S}_U$ be the set of all weighted unrooted trees without a specified limb $L$, let $S_U(x)$ be the generating function of $\mathcal{S}_U,$ and let $S_W$ be the generating function for weighted rooted trees without $L$ as a limb, then
$$S_U(x)=S_W(x)+\frac{1}{2}(S_W(x^2)-(S_W(x))^2).$$
Lu \cite{Lu1996} proved the analogous equation for the generation function of trees without a specified limbs and rooted tress without a specified limb. The proof for the weighted case is identical since the proof only concerns the structure of trees without considering their vertex weights.

Suppose the radius of convergence of a power series $A(x)$ is less than 1. Then the radius of convergence of $A(x^2)$ is less than $1$ but strictly greater than that of $A(x)$, and the radius of convergence of $(A(x))^2$ is the same as that of $A(x)$. Therefore, the radius of convergence of $S_W(x)+\frac{1}{2}(S_W(x^2)-(S_W(x))^2)$ is the same as that of $A(x)$. Consequently, the radius of convergence of $W(x)$ is $\alpha_T$ and the radius of convergence of $S_U(x)$ is $\alpha_S$. Theorem \ref{swtwrds} proved $\alpha_S>\alpha_T$. Similarly to Corollary \ref{haslimbwr}, we get the following result.

\begin{cor}
	For any given weighted tree $L$, the proportion of $n$-vertex weighted trees with $L$ as a limb approaches $1$ as $n$ approaches infinity.
\end{cor}

Therefore, we have extended Schwenk's result that almost all trees contain a specified limb to weighted trees.

\section{Application to quantum field theory}

Weighted rooted trees have many applications across combinatorics, mathematics, and beyond.  However, there is one particular application which has been inspirational for us and informs some of the choices we make.

In perturbative quantum field theory one wishes to understand particle interactions by looking at all the possible ways that some incoming particles could become some other set of outgoing particles, potentially with many intermediate particles.  The amplitude of this process can be calculated by summing the Feynman integrals for each of these possible ways.  Each of these ways can be represented by a graph, called the Feynman diagram or Feynman graph, and the graph determines the Feynman integral.  Unfortunately, these integrals are typically divergent, and physicists developed the process of renormalization to extract answers out of these divergent integrals.  These answers match experiment extremely well.

One way to mathematize renormalization is by using some combinatorial Hopf algebras known as renormalization Hopf algebras.  When the Feynman integral is already divergent upon integrating only the variables corresponding to a subgraph of the Feynman graph, then that subgraph is known as a subdivergence. The key thing that needs to be done for renormalization is to organize the structure of subdivergences within a Feynman graph.  This can be done by building a Hopf algebra on the graphs with a coproduct which pulls out subdivergences.  The original Hopf algebraic approach to renormalization, however, used rooted trees to represent the subdivergence structure of a Feynman diagram.  This gives the Connes-Kreimer Hopf algebra of rooted trees, see \cite{ck0}.  For references to further details see \cite{Yeats2017}.

Specifically, if any two subdivergences of a Feynman graph are either disjoint or contained one within the other then they directly have a rooted tree structure by containment.  Each vertex of the rooted tree then corresponds to a divergent subgraph, with the root corresponding to the whole graph.  Typically, each vertex of the tree is not labelled with this whole divergent subgraph, but rather with the result of taking this subgraph and contracting the subgraphs corresponding to its children in the tree.  This is called an insertion tree.  To recreate the graph from the insertion tree, one needs to carry around the additional information of where each graph should be inserted into its parent.

Subdivergences which are neither disjoint nor contained one within the other are known as overlapping.  Feynman graphs with overlapping subdivergences do not have unique insertion trees.  However, in each case with overlapping one can make a choice of one of the overlapping divergences and build an insertion tree using that subdivergence. For the purposes of renormalization, a Feynman graph can be replaced by the formal sum of its insertion trees built in this way.

For many purposes we don't need to carry around all the information in the insertion tree.  The most important piece of information to keep is the size of each graph.  The best measure of size for quantum field theory is the dimension of the cycle space of the graph, which is known as the loop number.  Consequently, the most important kind of trees are those which come from insertion trees but where the graphs at each vertex (the ones formed by contracting the children in the original subgraphs) are forgotten and only their loop number remains.  This results in a weighted rooted tree.  The loop number of the original Feynman diagram is the sum of the weights of the vertices of the tree, so summing the weights gives the correct notion of size for the tree.

Interpreting Schwenk's result in this context, a limb is a divergent subgraph (itself possibly with further subdivergences within it).  The fact that almost all weighted rooted tress have any given limb means that we cannot hope to avoid any particular subdivergence structure.  If some particular subdivergence is bad, then we are none-the-less stuck with it almost all the time.

The Connes-Kreimer Hopf algebra itself is defined as follows.  If $\mathcal{T}$ is the set of weighted rooted trees without an empty tree, then as an algebra the Connes-Kreimer Hopf algebra is $\mathbb{Q}[\mathcal{T}]$.  Writing $t_v$ for the subtree of $t$ rooted at $v\in V(t)$, then the coproduct on $t\in \mathcal{T}$ is
\[
\Delta(t) = \sum_{\substack{S\subset V(t)\\S\text{ antichain}}} \left(\prod_{v\in S}\right) t_v \otimes \left(t - \prod_{v\in S} t_v\right)
\]
with the convention that if $S$ is just the root then the right hand side of the tensor is interpreted to be $1$ rather than $0$.  The coproduct is extended to the whole algebra as an algebra homomorphism.  The counit is given by $\epsilon(1)=1$, $\epsilon(t)=0$ for $t\in \mathcal{T}$ and extended as an algebra homomorphism.  This gives a bialgebra, and it is graded by the (weighted) size of the tree with the $0$-graded piece simply being the underlying field.  Thus the antipode can be obtained recursively for free and we get a Hopf algebra, the Connes-Kreimer Hopf algebra.

The subtrees which can appear on the left hand side in the coproduct inspire our next definition of maximal limbs.

\section{Replacement of Maximal Limb}\label{sec max limb}

So far, the limbs we have considered are generic, in the sense that they are not chosen in any specific way. In this section, we consider a particular type of limbs named maximal limbs, and extend Schwenk's results to them. For a rooted tree $T$ and any vertex $v$ in $T$, the \emph{maximal limb} at $v$ is the limb at $v$ that contains exactly all descendants of $v$. This limb is maximal in a rooted setting because all limbs at $v$ are its subtrees. In this section, we prove a result that is similar to Theorem \ref{rtsamenum}, but for maximal limbs.

\begin{theorem}
	Let $L_1$ and $L_2$ be two distinct rooted trees with $\ell$ vertices, let $\ell_{1, n}$ and $\ell_{2, n}$ denote the number of $n$-vertex rooted trees without $L_1$ or $L_2$ as a maximal limb, respectively. Then $\ell_{1, n}=\ell_{2, n}$.
	\label{rtsamenummax}
\end{theorem}

\proof When $n<\ell$, the statement is trivially true because no rooted tree with $n$ vertices can have a limb with $l$ vertices. So we may assume $n\ge \ell$.

Use $\mathcal{T}_1$ to denote the set of $n$-vertex rooted trees with $L_2$ as a limb but without $L_1$ as limb, and $\mathcal{T}_2$ to denote the set of $n$-vertex rooted trees with $L_1$ as a limb but without $L_2$ as limb. It suffices to show that $|\mathcal{T}_1|=|\mathcal{T}_2|$. We prove this bijectively.

Let $T\in \mathcal{T}_1$, and consider the occurrences of the maximal limb $L_2$ in $T$. By definition of maximal limbs, these occurrences must be rooted at distinct vertices. If two occurrences of $L_2$ rooted at two distinct vertices $v_1$, $v_2$ with $i\ne j$ share a vertex $u$, then the union of the paths $v-v_1-u$ and $v-v_2-u$ contains a cycle, which is a contradiction. Therefore, no two of these occurrences of $L_2$ share a vertex. To obtain a tree in $\mathcal{T}_2$ from $t$, we simply replace each occurrence of the limb $L_2$ with the limb $L_1$. A similar argument can be applied to the trees in $\mathcal{T}_2$.

Therefore, this process is reversible, and we have constructed a bijection between $\mathcal{T}_1$ and $\mathcal{T}_2$, which completes the proof.\qed

Therefore, the number of $n$-vertex trees with a specified $\ell$-vertex limb is independent of the structure of the limb. In the next section, we show that almost all rooted trees contain a specified maximal limb by proving a number of more general statements.

\section{Maximal Limbs in Rooted Trees}

In the previous section, we extended Schwenk's limb replacement result to maximal limbs. Moreover, we can also ask, given a specified rooted tree $L$, would almost all rooted trees contain it as a maximal limb? In this section, we show that almost all trees contain a specified rooted tree as a limb.

Similar to the approach we used for weighted rooted trees, we start by deriving a recursion for the generating function of the set of rooted trees without $L$ as a maximal limb. For any $\ell$-vertex rooted tree $L$, let $\mathcal{S}$ be the set of rooted trees without $L$ as a maximal limb, and let $S(x)=\sum_{i\ge0}S_ix^i$ be the generating series of $\mathcal{S}$.

Given a rooted tree $T$ in $\mathcal{S}$. Delete the root of $T$, and observe that $F(T)$ is a set of rooted trees in $S_R$ that is not identical to $F(L)$.
Then $S(x)$ satisfies the following relation.

$$\mathcal{S}=(\mathcal{R}\times\mathbf{MSET}(\mathcal{S})) - (\mathcal{R}\times\mathcal{F(L)}).$$

The set subtraction is legitimate as $F(L) \in \mathbf{MSET}(\mathcal{S})$.

Therefore, the generating function of rooted trees without $L$ as a maximal limb satisfies the following recursive equation.


\begin{eqnarray*}
	S^*(x)&=&x\exp\left(\sum_{i=1}^\infty\frac{1}{i}S^*(x^i)\right)-x^\ell
\end{eqnarray*}

Using the same approach that we took in the case of weighted rooted trees, we denote the radius of convergence of $S(x)$ as $\alpha_S$, and we prove $S(\alpha_S)=1$.

\begin{lem}
	The series $S(x)$ satisfies $S(\alpha_S)+(\alpha_S)^\ell=1$.
	\label{radS3}
\end{lem}

\proof Define a multivariate function: for $x, y\in\cx$, $$F(x, y)=x\exp{\left(y+\sum_{i=2}^\infty\frac{1}{i}S(x^i)\right)}-x^\ell-y.$$

By the Implicit Function Theorem, $y=S(x)$ is the unique analytic solution of $F(x, y)=0$. Moreover, it has a singularity at $x=\alpha_S$, and $F(\alpha_S, S(\alpha_S))=0$.

Therefore, 
\begin{eqnarray*}
\frac{\partial F}{\partial y}(\alpha_S, S(\alpha_S))&=&F(\alpha_S, S(\alpha_S))+(\alpha_S)^\ell+S(\alpha_S)-1\\
&=&(\alpha_S)^\ell+S(\alpha_S)-1=0,
\end{eqnarray*}
due to the singularity. Thus $S(\alpha_S)+(\alpha_S)^\ell=1$. \qed

Meanwhile, recall that the generating function for rooted trees satisfies $$T(x)=x\exp\left(\sum_{i=1}^\infty\frac{1}{i}T(x^i)\right).$$ Let $r_T$ be the radius of convergence of $T(x)$.

\begin{theorem}
	The radius of convergence of $S(x)$ is greater than the radius of convergence of $T(x)$.
\end{theorem}

\proof It is known that $T(r_T)=1$ \cite{HP1973}, while Lemma \ref{radS3} gives $S(\alpha_S)+(\alpha_S)^\ell=1$.

By definition, $S(x)$ is coefficient-wise bounded above by $T(x)$, which implies $\alpha_S\ge r_T$. Moreover, since $T(x)$ enumerates all rooted trees with $L$ as a maximal limb and $S(x)$ does not, $\langle x^i, S(x)+x^\ell\rangle=\langle x^i, T(x)\rangle$ for all $i\le \ell$, and $\langle x^i, S(x)+x^\ell\rangle<\langle x^i, T(x)\rangle$ for all $i>\ell$. Therefore, $S(r_T)+(r_T)^\ell<T(r_T)=1$, which means $r_S\ne r_T$. Consequently, we obtain $\alpha_S>r_T$. \qed

Since the radius of convergence of $S(x)$ is greater than the radius of convergence of $T(x)$, the ratio between $\langle x^n, S(x)\rangle$ and $\langle x^n, T(x)\rangle$ approaches zero as $n$ gets arbitrarily large. Therefore, the proportion of $n$-vertex rooted trees with a specified maximal limb approaches $1$ as $n$ approaches infinity, so almost all rooted trees contains a specified maximal limb.

Note that the argument used in this section would also work for weighted rooted trees as defined in Section \ref{sec weighted limbs}. In other words, the proportion of weighted rooted trees with a specified maximal limb approaches $1$ as $n$ approaches infinity

\section{Constructing Graphs with Many Cospectral Vertices}







In this section, we introduce some preliminary results about cospectral vertices. In fact, there exist graphs with cospectral but not similar vertices\footnote{Two vertices are \emph{similar} if there is an automorphism of the graph taking one to the other.}. We will use the various characteristic polynomial identities to provide a step-by-step construction of trees with an arbitrarily large number of non-similar cospectral vertices, and to prove almost all trees have $k$ non-similar cospectral vertices for any integer $k\ge 2$.

\begin{figure}
	\centering
	\begin{tikzpicture}
	\node (1) at (0,0){};
	\node[label=$a$] (2) at (1,0){};
	\node (3) at (2,0){};
	\node (4) at (3,0){};
	\node[label=$b$] (5) at (4,0){};
	\node (6) at (5,0){};
	\node (7) at (6,0){};
	\node (8) at (7,0){};
	\node (9) at (2,1){};
	
	\path   (1) edge (2)
	(2) edge (3)
	(3) edge (4)
	(4) edge (5)
	(5) edge (6)
	(6) edge (7)
	(7) edge (8)
	(9) edge (3);
	\end{tikzpicture}
	\caption{The smallest tree with a pair of non-similar cospectral vertices $a$ and $b$}
	\label{schwenkstree}
\end{figure} 

Figure \ref{schwenkstree} is the smallest tree with a pair of non-similar cospectral vertices, $a$ and $b$, as labeled in the figure. Note that this is the tree Schwenk used to perform his limb replacement proof that almost all trees have a cospectral mate \cite{Schwenk1973}. The pair of non-similar cospectral vertices are exactly the pair of vertices he chose as the roots of his limbs.  Any tree with a pair of non-similar cospectral vertices would work for Schwenk's proof that almost all trees have a cospectral mate, where the cospectral vertices serve as roots of the limbs; in fact any time there is a set of pair-wise cospectral vertices in the tree $L$, we can use different rootings of $L$ as a limb to build cospectral trees. This is a consequence of Lemma~\ref{onesumid}.

To further discuss cospectral vertices, the following identities from Godsil's book are useful 
\cite[Theorem 2.1.5]{Godsil1993}. 
\begin{lem}
	\begin{enumerate}[(a)]
		\item $\phi(G\cup H)=\phi(G)\phi(H)$.
		\item $\phi(G)=\phi(G\backslash e)-\phi(G\backslash\{u, v\})$ if $e=uv$ is a cut-edge of $G$.
		\item $\frac{d}{dx}\phi(G)=\sum_{i\in V(G)}\phi(G\backslash i)$.
	\end{enumerate}
	\label{chariden}
\end{lem}

Note that $G\backslash e$ denotes the graph obtained by only deleting the edge $e$, without deleting the end vertices of $e$. On the other hand, $G\backslash\{u, v\}$ denotes the graph obtained by deleting the vertices $u$ and $v$, as well as all edges incident to $u$ or $v$. These two notations are used throughout our discussions in our paper.


Any tree $T$ that contains $L$ as a limb can be constructed by identifying a vertex in the limb $L$ with a vertex in another tree. We call the identification operation $1$-sum. Formally, let $H, K$ be two graphs and $u\in H$, $v\in K$. Then the \emph{$1$-sum} of $H$ and $K$, denoted as $H\onesum K$, is the graph obtained by identifying a vertex in $H$ and a vertex in $K$. 

Note that in Schwenk's original proof to Theorem \ref{treecosmate}, he used a limb replacement technique, which is in fact a $1$-sum. He showed that there exist two trees $R$ and $S$ and two specific vertices $u\in R$ and $v\in S$, such that the $1$-sum of any tree $T$ at any vertex $w\in T$ with $R$ at $u$ is cospectral with the $1$-sum of $T$ at $w$ with $R$ at $u$. In general, we have the following $1$-sum identity for the characteristic polynomial \cite[Corollary 4.3.3]{Godsil1993}.

\begin{lem}
	Let $G=H\onesum K$, with the identified vertex being $v$. Then $$\phi(G)=\phi(H\backslash v)\phi(K)+\phi(H)\phi(K\backslash v)-x\phi(H\backslash v)\phi(K\backslash v).$$
	\label{onesumid}
\end{lem}

\vspace{-0.5cm}

The identities in Lemma \ref{chariden} and Lemma \ref{onesumid} play essential roles in our construction of trees with cospectral vertices below. Moreover, they will be extended to weighted graphs and weighted rooted trees in the next section. 

To construct trees with many cospectral vertices, we start by introducing a graph operation that we named one-vertex extension. Given a graph $G$ and a set of vertices $S\subseteq V(G)$, we define the \emph{one-vertex extension} of $G$ with respect to $S$ to be the graph obtained by adding a vertex to $G$ and connecting this vertex to all vertices in $S$ and nothing else. Note that when $S=V(G)$, the one-vertex extension with respect to $S$ is simply the cone over $G$.

Let $F$ be a graph with exactly $c$ components, such that its components are pairwise cospectral. Suppose there exists distinct sets $S_1, S_2, \ldots, S_m \subset V(F)$ of pairwise cospectral vertices for some integer $m\ge2$. Moreover, suppose that for all $1\le i\le m$, we have $|S_i|=c$ and each $S_i$ contains exactly one vertex from each component. 

Observe that such a graph $F$ exists. For example, take two disjoint copies of the graph in Figure \ref{schwenkstree} as $F$, with the two isomorphic components being $C_1$, $C_2$. If we call the pair of cospectral vertices in $C_i$ as $a_i, b_i$ for $i=1, 2$. Let $S_1=\{a_1, a_2\}$, $S_2=\{b_1, a_2\}$, and $S_3=\{b_1, b_2\}$. This graph satisfies the condition listed in the paragraph above with $c=2$ and $m=3$. Moreover, in this example, the graphs $F\backslash S_i$ are pairwise non-isomorphic. This fact will be further discussed later in this section.

In general, given a graph $G$ and two vertices $u, v\in G$, if $G\backslash u\cong G\backslash v$ and there is no automorphism of $G$ mapping $u$ to $v$, then $u$ and $v$ are a pair of \emph{pseudo-similar} vertices. Any graph $G$ with a pair of pseudo-similar vertices can be used to construct $F$, where $F$ consists of two copies of $G$. In this case, $c=2$ and $m=3$ as well. Herndon and Ellezy gave a construction for graphs with pairs of pseudo-similar vertices in 1975 \cite{HE1975}. Later, Godsil and Kocay proved all graphs with a pair of pseudo-similar vertices can be constructed by Herndon and Ellezy's construction \cite{GK1982}. Therefore, other than our example in the previous paragraph, there exists a construction for even more graphs satisfying the conditions of $F$, with $c=2$ and $m=3$. Later in this section, we show our construction can create $F$ with greater values of $c$ and $m$.

The following result is used in the main proof.

\begin{theorem}
	Let $G$ be a graph with a set of cospectral vertices $A=\{a_1, a_2, \cdots, a_k\}$, where each component of $G$ contains at most one vertex in $A$. Let $S_1, S_2\subseteq A$ such that $|S_1|=|S_2|$. For $i=1, 2$, let $G_i$ be the one-vertex extension of $G$ with respect to $S_1$ and $S_2$. Then $G_1$ and $G_2$ are cospectral.
	\label{cmgphsfromcmvtcs}
\end{theorem}

\proof We prove this theorem by induction on the size of $S_1$ and $S_2$. Suppose $|S_1|=|S_2|=1$. Without loss of generality, let $v$ be connected to $a_i$ in $G_i$, for $i=1, 2$. Use $e_{a_iv}$ to denote the edge with $a_i$ and $v$ being its endpoints. Each $e_{a_iv}$ is a cut edge and so 

\begin{align*}
	\phi(G_1)&=\phi(G_1\backslash e_{a_1v})-\phi(G_1\backslash\{a_1, v\})\\
	&=x\phi(G)-\phi(G\backslash a_1)\\
	&=\phi(G_2\backslash e_{a_2v})-\phi(G\backslash a_2)\\
	&=\phi(G_2).
\end{align*}

Now assume the theorem statement holds for all $|S_1|=|S_2|\le \ell-1$, where $1<\ell\le k$. Without loss of generality, suppose $|S_1|=|S_2|=k$, $a_1\in S_1$, $a_2\in S_2$. Since each component of $G$ contains at most one vertex in $A$, $e_{a_iv}$ is a cut edge. By the induction hypothesis, $\phi(G_1\backslash e_{a_1v})=\phi(G_2\backslash e_{a_2v})$.

Note that 
\begin{align*}
	\phi(G_1)&=\phi(G_1\backslash e_{a_1v})-\phi(G_1\backslash\{a_1, v\})\\
	&=\phi(G_1\backslash e_{a_1v})-\phi(G\backslash a_1)\\
	&=\phi(G_2\backslash e_{a_2v})-\phi(G\backslash a_2)\\
	&=\phi(G_2).
\end{align*}

Hence the theorem statement is true. \qed

\begin{figure}
	\centering
	\begin{minipage}{0.45\textwidth}
		\begin{tikzpicture}[xscale=0.38, yscale=0.5]
			\node[circle,fill=black] (1) at (0,0){};
			\node[circle,fill=black, label=$a_1$] (2) at (-1,1){};
			\node[circle,fill=black] (3) at (-2,0){};
			\node[circle,fill=black] (4) at (-3,-1){};
			\node[circle,fill=black] (5) at (-4,-2){};
			\node[circle,fill=black] (6) at (-5,-3){};
			\node[circle,fill=black] (7) at (-6,-4){};
			\node[circle,fill=black] (8) at (-7,-5){};
			\node[circle,fill=black] (9) at (-1,-1){};
			
			\path   (1) edge (2)
			(2) edge (3)
			(3) edge (4)
			(4) edge (5)
			(5) edge (6)
			(6) edge (7)
			(7) edge (8)
			(9) edge (3);
			
			\node[circle,fill=black] (11) at (8,-3){};
			\node[circle,fill=black] (12) at (7,-2){};
			\node[circle,fill=black] (13) at (6,-1){};
			\node[circle,fill=black] (14) at (5,0){};
			\node[circle,fill=black, label=$b_2$] (15) at (4,1){};
			\node[circle,fill=black] (16) at (3,0){};
			\node[circle,fill=black] (17) at (2,-1){};
			\node[circle,fill=black] (18) at (1,-2){};
			\node[circle,fill=black] (19) at (5,-2){};
			
			\path   (11) edge (12)
			(12) edge (13)
			(13) edge (14)
			(14) edge (15)
			(15) edge (16)
			(16) edge (17)
			(17) edge (18)
			(19) edge (13);
			
			\node[circle,fill=black] (10) at (1.5,3){};
			
			\path   (10) edge (2)
			(15) edge (10);
		\end{tikzpicture}
	\end{minipage}
	\qquad
	\begin{minipage}{0.45\textwidth}
		\begin{tikzpicture}[xscale=0.38, yscale=0.5]
			\node[circle,fill=black] (1) at (0,0){};
			\node[circle,fill=black, label=$a_1$] (2) at (-1,1){};
			\node[circle,fill=black] (3) at (-2,0){};
			\node[circle,fill=black] (4) at (-3,-1){};
			\node[circle,fill=black] (5) at (-4,-2){};
			\node[circle,fill=black] (6) at (-5,-3){};
			\node[circle,fill=black] (7) at (-6,-4){};
			\node[circle,fill=black] (8) at (-7,-5){};
			\node[circle,fill=black] (9) at (-1,-1){};
			
			\path   (1) edge (2)
			(2) edge (3)
			(3) edge (4)
			(4) edge (5)
			(5) edge (6)
			(6) edge (7)
			(7) edge (8)
			(9) edge (3);
			
			\node[circle,fill=black] (11) at (3,0){};
			\node[circle,fill=black, label=$a_2$] (12) at (4,1){};
			\node[circle,fill=black] (13) at (5,0){};
			\node[circle,fill=black] (14) at (6,-1){};
			\node[circle,fill=black] (15) at (7,-2){};
			\node[circle,fill=black] (16) at (8,-3){};
			\node[circle,fill=black] (17) at (9,-4){};
			\node[circle,fill=black] (18) at (10,-5){};
			\node[circle,fill=black] (19) at (4,-1){};
			
			\path   (11) edge (12)
			(12) edge (13)
			(13) edge (14)
			(14) edge (15)
			(15) edge (16)
			(16) edge (17)
			(17) edge (18)
			(19) edge (13);
			
			\node[circle,fill=black] (10) at (1.5,3){};
			
			\path   (10) edge (2)
			(12) edge (10);
		\end{tikzpicture}
	\end{minipage}
	\caption{A pair of cospectral graphs constructed using Theorem \ref{cmgphsfromcmvtcs}.}
	\label{csgphsfromcsvtcs}
\end{figure}

Figure \ref{csgphsfromcsvtcs} gives an example of a pair of cospectral graphs constructed using Theorem \ref{cmgphsfromcmvtcs}. Each of the two graphs in Figure \ref{csgphsfromcsvtcs} is a one-vertex extension of two copies of the tree in \ref{schwenkstree}, where the extension is done on different but cospectral vertices.

Note that for a vertex transitive graph $G$, all vertices of $G$ can be put in the set $A$ defined in Theorem \ref{cmgphsfromcmvtcs}. Therefore, this result gives a construction of an infinite number of large sets containing pairwise cospectral graphs.

\begin{lem}
	Let $F$ be a graph with $c$ components that are pairwise cospectral. For each $j=1, 2$, suppose there exists distinct sets $S^j_1, S^j_2, \ldots, S^j_m \subseteq V(F)$ of pairwise cospectral vertices for some integer $m\ge1$. Moreover, for all $1\le i\le m$ and $j=1, 2$, we have $|S^j_i|=c$, and each $S^j_i$ contains exactly one vertex from each component. 
	
	For $1\le i\le m$ and $j=1, 2$, let $F^j_i$ be the one-vertex extension of $F$ with respect to $S^j_i$, with the added vertex being $v^j_i$. Let $H_j$ be the one-vertex extension of $\cup_{i=1}^mF^j_i$ with respect to $v^j_1, \ldots, v^j_m$, with the added vertex being $r_j$. Then $\phi(H_1)=\phi(H_2)$.
	\label{cmtrees1sumk1}
\end{lem}

\proof We prove this result by induction on $m$.

When $m=1$, let $e_1$ be the only edge $r_1$ is incident to, and let $e_2$ be the only edge $r_2$ is incident to. Then $v^1_1$ and $v^2_1$ are the other end-vertices of $e_1$ and $e_2$, respectively. Note that by Theorem \ref{cmgphsfromcmvtcs}, $$\phi(H_1\backslash e_1)=\phi(H_2\backslash e_2),$$ and by construction, $$\phi(H_1\backslash\{r_1, v^1_1\})=\phi(H_2\backslash\{r_2, v^2_1\}).$$ Lemma \ref{chariden} (b) gives $\phi(H_1)=\phi(H_2)$.

Now suppose the lemma holds when $m=n-1$ some $n>=1$. We consider the case where $m=n$. Let $e_1$ and $e_2$ be edges that connect $r_1$ and $r_2$ to $v^1_1$ and $v^2_1$, respectively. Note that by the induction hypothesis and Theorem \ref{cmgphsfromcmvtcs}, $$\phi(H_1\backslash e_1)=\phi(H_2\backslash e_2),$$ and by the construction, $$\phi(H_1\backslash\{r_1, v^1_1\})=\phi(H_2\backslash\{r_2, v^2_1\}).$$ Then once again Lemma \ref{chariden} (b) gives $\phi(H_1)=\phi(H_2)$. \qed

Now we have all the tools to prove the following main result of this section. We show that, by applying the one-vertex extension twice in a row on a graph satisfying some conditions, a graph with cospectral vertices is constructed.

\begin{theorem}
	Let $F$ be a graph with exactly $c$ components that are pairwise cospectral. Suppose there exists distinct sets $S_1, S_2, \ldots, S_m \subseteq V(F)$ of pairwise cospectral vertices for some integer $m\ge2$. Moreover, for all $1\le i\le m$, we have $|S_i|=c$, and each $S_i$ contains exactly one vertex from each component. 
	
	For $1\le i\le m$, let $F_i$ be the one-vertex extension of $F$ with respect to $S_i$, with the added vertex being $v_i$. Let $G$ be the one-vertex extension of $\cup_{i=i}^mF_i$ with respect to $v_1, \ldots, v_m$, with the added vertex being $v$. Then $\cup_{i=1}^mN(v_m)\backslash v$ is a set of pairwise cospectral vertices in both $G$ and $G\backslash v$..
	
	Moreover, if the graphs $F\backslash S_i$ are pairwise non-isomorphic for $1\le i\le m$, then the vertices in the set $\cup_{i=1}^mN(v_m)\backslash v$ are pairwise non-similar in $G$.
	\label{csvtc_ext_thm}
\end{theorem}

\proof Let $u_1, u_2\in \cup_{i=1}^mN(v_m)\backslash v$. When the graphs $F\backslash S_i$ are pairwise non-isomorphic, $u_1$ and $u_2$ are not similar. 
We first show that $\phi(G\backslash u_1)=\phi(G\backslash u_2)$.

Let $1\le i_1, i_2\le m$ such that $v_{i_1}$ is a neighbor of $u_1$ and $v_{i_2}$ is a neighbor of $u_2$. Note it is possible for $i_1=i_2$. Use $e_1$, $e_2$ to denote the edges between $v_{i_1}$ and $v$ or $v_{i_2}$ and $v$, respectively. Let $a_1, a_2$ be vertices in $F$ such that there is an isomorphism taking $a_j$ in $F$ to $u_j$ in $F_{i_j}\backslash v_{i_j}$, for $j=1, 2$.

We first compare $\phi(G\backslash\{u_1, v_{i_1}, v\})$ to $\phi(G\backslash\{u_2, v_{i_2}, v\})$. For $j=1, 2$, observe that $$G\backslash\{u_j, v_{i_j}, v\}=(\left(\cup_{\ell=1}^mF_\ell\right)\backslash F_{i_j}) \cup (F\backslash a_j).$$

Since $u_1$ and $u_2$ originally came from cospectral vertices in $F$, $\phi(F\backslash a_1)=\phi(F\backslash a_2)$. Moreover, by Theorem \ref{cmgphsfromcmvtcs}, the graphs $F_1, \ldots, F_m$ are pairwise cospectral. Therefore, $$\phi(G\backslash\{u_1, v_{i_1}, v\})=\phi(G\backslash\{u_2, v_{i_2}, v\}).$$

Now we compare $\phi((G\backslash u_1)\backslash e_1)$ with $\phi((G\backslash u_2)\backslash e_2)$. For $j=1, 2$, use $B_j$ to represent the copy of $F$ in $G$ containing the vertex $u_j$, then $$(G\backslash u_i)\backslash e_i=(G\backslash F_{i_j})\cup(F_{i_j}\backslash B_j)\cup(F\backslash a_j).$$

As discussed above, $\phi(F\backslash a_1)=\phi(F\backslash a_2)$. Meanwhile, Theorem \ref{cmgphsfromcmvtcs} implies $\phi(F_{i_1}\backslash B_1)=\phi(F_{i_2}\backslash B_2)$. Then $\phi(G\backslash F_{i_1})=\phi(G\backslash F_{i_2})$ is a result of Lemma \ref{cmtrees1sumk1}. Therefore, $$\phi((G\backslash u_1)\backslash e_1)=\phi((G\backslash u_2)\backslash e_2)$$.

By Lemma \ref{chariden} (b), $\phi(G\backslash u_1)=\phi(G\backslash u_2)$.

To see $u_1$ and $u_2$ are cospectral vertices in $G\backslash v$, observe that, for $j=1, 2$, we have $$G\backslash\{u_j, v\}=(\left(\cup_{\ell=1}^mF_\ell\right)\backslash F_{i_j}) \cup (F_{i_j}\backslash u_j)=(\left(\cup_{\ell=1}^mF_\ell\right)\backslash F_{i_j}) \cup (F_{i_j}\backslash B_j) \cup (F\backslash a_j).$$ We've already proven $\phi((\left(\cup_{\ell=1}^mF_\ell\right)\backslash F_{i_1}))=\phi((\left(\cup_{\ell=1}^mF_\ell\right)\backslash F_{i_2}))$, $\phi(F_{i_1}\backslash B_1)=\phi(F_{i_2}\backslash B_2)$, and $\phi(F\backslash a_1)=\phi(F\backslash a_2)$. Therefore $u_1$ and $u_2$ are cospectral vertices in $G\backslash v$. \qed

\begin{figure}
	\centering
	\begin{tikzpicture}
		\draw (0,0) ellipse (20pt and 35pt) node[fill=white] {$L$};
		\draw (2,0) ellipse (20pt and 35pt) node[fill=white] {$L$};
		\draw (4,0) ellipse (20pt and 35pt) node[fill=white] {$L$};
		\draw (6,0) ellipse (20pt and 35pt) node[fill=white] {$L$};
		\draw (8,0) ellipse (20pt and 35pt) node[fill=white] {$L$};
		\draw (10,0) ellipse (20pt and 35pt) node[fill=white] {$L$};
		
		\node[circle,fill=black] (1) at (0,0.5){};
		\node[circle,fill=black] (2) at (2,0.5){};
		\node[circle,fill=black] (3) at (4,0.5){};
		\node[circle,fill=black] (4) at (6,0.5){};
		\node[circle,fill=black] (5) at (8,0.5){};
		\node[circle,fill=black] (6) at (10,0.5){};
		
		\node[circle,fill=white, scale=0.05, label=$a_1$] (7) at (-0.2,0.5){};
		\node[circle,fill=white, scale=0.05, label=$a_2$] (8) at (2.2,0.5){};
		\node[circle,fill=white, scale=0.05, label=$a_1$] (9) at (3.8,0.5){};
		\node[circle,fill=white, scale=0.05, label=$b_2$] (10) at (6.2,0.5){};
		\node[circle,fill=white, scale=0.05, label=$b_1$] (11) at (7.8,0.5){};
		\node[circle,fill=white, scale=0.05, label=$b_2$] (12) at (10.2,0.5){};
		
		\node[circle,fill=black] (13) at (1,2){};
		\node[circle,fill=black] (14) at (5,2){};
		\node[circle,fill=black] (15) at (9,2){};
		
		\draw (1)--(13);
		\draw (2)--(13);
		\draw (3)--(14);
		\draw (4)--(14);
		\draw (5)--(15);
		\draw (6)--(15);
		
		\node[circle,fill=black] (16) at (5,4){};
		\draw (16)--(13);
		\draw (16)--(14);
		\draw (16)--(15);
	\end{tikzpicture}
	\label{exconst}
	\caption{A graph obtained by applying our construction in this section to Figure \ref{schwenkstree}.}
\end{figure}

As shown in Figure \ref{schwenkstree}, there exists a connected graph with a pair of non-similar cospectral vertices. Let this graph be $L$. Furthermore, as discussed earlier in this section, by taking $c=2$, we can get three sets of two vertices $S_1$, $S_2$, and $S_3$ such that $F\backslash S_1$, $F\backslash S_2$, and $F\backslash S_3$ are pairwise non-isomorphic. Therefore, our extension would generate a connected graph as shown in Figure \ref{exconst}, which has six cospectral vertices, four of which are non-similar. In general, we can treat the sets $S_i$ as distinct multisets of cospectral vertices in the graphs $F$. Then if there exists a connected graph $G$ with $k$ non-similar cospectral vertices, for an arbitrary positive integer value $c\ge2$ we can construct a connected graph with $c\binom{k+c-1}{c}$ cospectral vertices, where at least $\binom{k+c-1}{c}$ of them are non-similar.

If we apply this construction repeatedly, graphs with an arbitrarily large number of non-similar cospectral vertices can be constructed. Moreover, since this construction does not create cycles in the graph, and there exists a tree (Schwenk's tree) with a pair of non-similar cospectral vertices, we can say that for any integer $k\ge2$, there exists a tree with $k$ cospectral vertices.

This construction can be applied to graphs that are not trees as well. For example, Figure \ref{cospvtcs} gives a graph that is not a tree and has a pair of non-similar cospectral vertices.

\begin{cor}
	For any $k\geq 2$ there exist trees with $k$ cospectral vertices.
\end{cor}

On the other hand, there exists graphs that are not trees with a pair of non-similar cospectral vertices. An example is shown in Figure \ref{cospvtcs}, where $a$ and $b$ are a pair of non-similar cospectral vertices. Consequently, we get a similar result for that are not trees.

\begin{figure}
	\centering
	\begin{tikzpicture}
		\node (1) at (0,0){};
		\node (2) at (1,0){};
		\node (3) at (2,0){};
		\node (4) at (3,0){};
		\node (5) at (4,0){};
		\node (6) at (5,0){};
		\node (7) at (0.5,0.75){};
		\node (8) at (2.5,0.75){};
		\node[circle, fill=white, scale=0.25, label=$a$] (9) at (1,-0.6){};
		\node[circle, fill=white, scale=0.25, label=$b$] (10) at (3,-0.6){};
		
		\path   (1) edge (2)
		(2) edge (3)
		(3) edge (4)
		(4) edge (5)
		(5) edge (6)
		(1) edge (7)
		(7) edge (2)
		(8) edge (3)
		(8) edge (4);
	\end{tikzpicture}
	\caption{A graph with cycles and a pair of non-similar cospectral vertices $a$ and $b$}
	\label{cospvtcs}
\end{figure} 

\begin{cor}
	For any $k\geq 2$ there exist graphs which are not trees with $k$ cospectral vertices.
\end{cor}

Furthermore, we can also obtain non-trees by additional 1-sums as described in the following result.  However, for both ways of obtaining non-trees, this construction generates rather special graphs, in particular with many cut vertices, and so is insufficient for a general picture about cospectral vertices in graphs.

\begin{theorem}
	Let $G$ be a graph. For any integer $k\ge 2$, suppose $L$ is a graph with $k$ cospectral vertices $a_1$, $a_2$, $\ldots$, $a_k$. Moreover, assume there exists a vertex $r\in L$ such that $r\not\in\{a_1, a_2, \ldots, a_k\}$, and the vertices $a_1$, $a_2$, $\ldots$, $a_k$ are still cospectral vertices in $L\backslash r$. Let $G'$ be the $1$-sum of $G$ and $L$ by identifying an arbitrary vertex in $G$ and $r$.  Then the vertices $a_1$, $a_2$, $\ldots$, $a_k$ are cospectral $G'$.
	\label{many_cos_vtc}
\end{theorem}

\proof For any $1\le i\le k$, by Lemma \ref{onesumid}, we get
\[\phi(G'\backslash a_i)=\phi(G)\phi(L\backslash \{a_i, r\})+\phi(G\backslash v)\phi(L\backslash a_i)-x\phi(G\backslash v)\phi(L\backslash \{a_i, r\}).\]

Then for any $1\le j\le k$, since $$\phi(L\backslash a_i)=\phi(L\backslash a_j)$$ and $$\phi(L\backslash \{a_i, r\})=\phi(L\backslash \{a_j, r\}),$$ we get $$\phi(G'\backslash a_i)=\phi(G'\backslash a_j).$$ In other words, $a_1, a_2, \ldots, a_k$ are cospectral vertices in $G'$.\qed

Note that by Theorem \ref{csvtc_ext_thm}, such a graph $L$ exists, and further such graphs $L$ exist which are trees and which are not trees.  In the case when both $G$ and $L$ are trees, taking the $1$-sum can be viewed as attaching $L$ as a limb of $G$ at the vertex $r$. Schwenk proved almost all trees have $L$ as a limb \cite{Schwenk1973}, which implies almost all trees have $k$ pairwise non-similar cospectral vertices.  This gives the following corollary. 

\begin{cor}
	For any $k\geq 2$, the proportion of $n$-vertex trees with $k$ cospectral vertices approaches $1$ as $n$ approaches infinity.
\end{cor}

However, even with more general input, the construction still results in a cut-vertex, while most graphs do not have cut-vertex.

In Section \ref{cos_weighted}, we proved our definition of weighted characteristic polynomial preserved a lot of properties of the characteristic polynomial for unweighted graphs. Consequently, our construction for cospectral vertices in graphs applies to weighted graphs as well: given weighted trees or graphs as input, satisfying the hypotheses as before but with weighted cospectrality in place of cospectrality, each time a 1-vertex extension is used in the construction simply assign weight $1$ to the new vertex. In particular then we can build weighted graphs $L$ with $k$ weighted cospectral vertices satisfying the conditions for Theorem \ref{many_cos_vtc}.  Then let $G$ be a weighed graph and $v\in V(G)$ with any positive integer weight. By Lemma \ref{onesumwr}, the $k$ weighted cospectral vertices in $L$ are still weighted cospectral in a weighted graph obtained by identifying $v$ in $G$ and $r$ in $L$, regardless of the specific weight of the identified vertex. Specifically, our results in this Section apply to weighted trees as well.

\begin{cor}
	For any $k\geq 2$, the proportion of $n$-vertex weighted trees with $k$ weighted cospectral vertices approaches $1$ as $n$ approaches infinity.
\end{cor}

\section{Cospectrality of Weighted Graphs}
\label{cos_weighted}

In this section, we prove that almost all weighted rooted trees are weighted cospectral with another weighted rooted tree.

We start by proving the weighted rooted trees satisfy identities that are analogous to Theorem \ref{chariden} and Lemma \ref{onesumid}. To do so, we need the following two results, where the proof of the former is in the book by Godsil \cite[Theorem 2.1.1]{Godsil1993}.

\begin{theorem}
	Let $X$ and $Y$ be any $n\times n$ matrices. Then $\det(X+Y)$ is equal to the sum of the determinants of the $2^n$ matrices obtained by replacing each subset of the columns of $X$ by the corresponding subset of the columns of $Y$.
	\label{detfor}
\end{theorem}

\begin{lem}
	Let $M$ be a $n\times n$ block matrix of the form
	\[
	\begin{bmatrix} C& \mathbf{0} \\ \mathbf{0} & D\end{bmatrix},
	\]
	where $C$ is a $\ell\times m$ matrix, $D$ is a $(n-\ell)\times (n-m)$ matrix, with $\ell\neq m$, and the matrices $\mathbf{0}$ denote all-zero matrices of appropriate dimensions. Then $\det(M)=0$.
	\label{mtxblocktri0det}
\end{lem}

\proof Without loss of generality, we may assume that $\ell>m$. Suppose $D=[D_1|D_2]$ where $D_1$ is a $(n-\ell)\times(\ell-m)$ matrix and $D_2$ is a $(n-\ell)\times(n-\ell)$ matrix. Then we can view $M$ as a block lower triangular matrix
	\[
\begin{bmatrix} C'& \mathbf{0} \\ E & D_2\end{bmatrix},
\]
where $C'=[C|\mathbf{0}]$ is a $\ell\times\ell$ matrix, $E=[\mathbf{0}|D_1]$ is a $(n-\ell)\times\ell$ matrix, with $\mathbf{0}$ again representing all-zero matrices of appropriate dimensions. Clearly, $\det(C')=0$. Then $\det(M)=\det(C')\det(D_2)=0$. \qed

Now we have all the necessary tools to prove the following result, analogous to the unweighted case. Recall that $G\backslash e$ denotes the graph obtained by only deleting the edge $e$, while $G\backslash\{u, v\}$ denotes the graph obtained by deleting the vertices $u$ and $v$, as well as all edges incident to $u$ or $v$.

\begin{lem}
	Let $W, W_1, W_2$ be weighted graphs. Then,
	\begin{enumerate}[(a)]
		\item $\phi^*(W_1\cup W_2)=\phi^*(W_1)\phi^*(W_2)$;
		\item $\phi^*(W)=\phi^*(W\backslash e)-\phi^*(W\backslash\{u, v\})$ if $e=uv$ is a cut-edge;
		\item $\frac{d}{dx}\phi^*(W)=\sum_{i\in V(W)}w(i)x^{w(i)-1}\phi^*(W\backslash i)$.
	\end{enumerate}
	\label{wtid}
\end{lem}

\proof
\begin{enumerate}[(a)]
	\item For square matrices $A$ and $B$, note that $\det\left[\begin{matrix}
	A&\mathbf{0}\\\mathbf{0}&B
	\end{matrix}\right]=\det(A)\det(B)$. Therefore, $$\phi^*(W_1\cup W_2)=\phi^*(W_1)\phi^*(W_2).$$
	\item Suppose $W$ has $n$ vertices. Let $E_{uv}$ be the $n\times n$ matrix whose all entries are zeros except that the $uv$-entry and $vu$-entry are $1$'s. Recall that $e$ is a cut-edge, so deleting it would result in two components. We may assume the $n\times n$ matrices in our discussion are arranged in a way where the first $\ell$ rows and columns correspond to the component $C_1$ with $u$ in it, with the $\ell$-th row and column corresponding $u$, and the $(\ell+1)$-st row and column corresponding to $v$. The other component is $C_2$. Observe that \[\phi^*(W)=\det(I^*-A(W))=\det(I^*-A+E_{uv}-E_{uv}).\] Let $X=I^*-A(W)+E_{uv}$ and $Y=-E_{uv}$, then we can apply Theorem \ref{detfor} to compute $\det(W)$. The matrix $X=I^*-A(W)+E_{uv}$ is illustrated in Figure \ref{i-a+e}.
	
	\begin{figure}
		\centering
		\[\left[\begin{array}{c|c|c|c}
		I_1^* - A(C_1\backslash u) & \cdots & \mathbf{0} & \mathbf{0}\\
		\hline
		\vdots & x^{w(u)} & 0 & \mathbf{0}\\
		\hline
		\mathbf{0} & 0 & x^{w(v)} & \vdots \\
		\hline
		\mathbf{0} & \mathbf{0} & \cdots & I_1^* - A(C_2\backslash v)
		\end{array}\right]\]
		\caption{The matrix $X=I^*-A(G)+E_{uv}$. The matrices $I_1^*$ and $I_2^*$ are diagonal matrices of appropriate dimensions, whose diagonal entries are $x^{w(i)}$, such that $i$ is the vertex corresponding to the respective row and column. The bolded $\mathbf{0}$'s represent the all-zero matrix of appropriate dimension.}
		\label{i-a+e}
\end{figure}

\begin{figure}		
		\[\left[\begin{array}{c|c|c|c}
		I^* - A(C_1\backslash u) & 0 & \mathbf{0} & \mathbf{0}\\
		\hline
		\vdots & 0 & 0 & \mathbf{0}\\
		\hline
		\mathbf{0} & -1 & x^{w(v)} & \vdots \\
		\hline
		\mathbf{0} & 0 & \cdots & I^* - A(C_2\backslash v)
		\end{array}\right]\]
		\caption{The matrix $B$. The matrices $I_1^*$ and $I_2^*$ are diagonal matrices of appropriate dimensions, whose diagonal entries are $x^{w(i)}$, such that $i$ is the vertex corresponding to the respective row and column. The bolded $\mathbf{0}$'s represent the all-zero matrix of appropriate dimension.}
		\label{matrixb}
\end{figure}

	Since any matrix with an all zero column has determinant 0, when using Theorem \ref{detfor} to replace columns of $I^*-A(W)+E_{uv}$, only three cases are considered: not replacing a column at all, replacing column $u$ and column $v$, and replacing exactly one of column $u$ and column $v$. For the first case, the determinant of the resulting matrix is $\det(I^*-A(W)+E_{uv})$, which is $\phi^*(W\backslash e)$. For the second case, by Laplace's formula (cofactor expansion) along the $\ell$-th and $(\ell+1)$-st columns we see that the determinant of the resulting matrix is $-\phi^*(W\backslash\{u, v\})$.
	
	Now we consider the last case. Assume first that the $\ell$-th column in $I^*-A(W)+E_{uv}$ is replaced by the $\ell$-th column of $-E_{uv}$. We call the resulting matrix is $B$, and we consider the determinant of $B$. The matrix $B$ is illustrated in Figure \ref{matrixb}.
	
        Let $B'$ be the matrix obtained from $B$ by deleting its $(\ell+1)$-st row and $\ell$-th column, then by cofactor expansion along the $\ell$-th column of $B$ we get $\det(B)=\det(B')$. Note that the entries in the $\ell$-th column of $B'$ come from the $(\ell+1)$-st column of $B$. In particular, the first $\ell$ entries in this column are zeros, and all other entries in this column can be $0$ or $-1$.  Therefore $B'$ is a block matrix of the form
        \[
        \begin{bmatrix} C& \mathbf{0} \\ \mathbf{0} & D\end{bmatrix}
        \] 
        where $C$ is $\ell\times (\ell-1)$ and $D$ is $(n-\ell-1)\times (n-\ell)$. By Lemma \ref{mtxblocktri0det}, $\det(B')=0$. An analogous argument applies when replacing only the $(\ell+1)$-st column. Therefore, the last case simply gives a determinant of $0$.
        
	
	
	Summing the determinants from all three cases, we get $$\phi^*(W)=\phi^*(W\backslash e)-\phi^*(W\backslash\{u, v\}).$$
	
	\item Let $S_n$ be the set of all permutations of $\{1, 2, \ldots, n\}$. Recall that for a matrix $M$, if $m_{ij}$ is its $ij$-entry, then its determinant can be computed by the Leibniz formula: \[\det(M)=\sum_{\sigma\in S_n}\left(\text{sgn}(\sigma)\prod_{i=1}^nm_{i, \sigma(i)}\right),\] where $\text{sgn}(\sigma)$ is the sign of the permutation. Use $c_{ij}$ to denote the $ij$-entry of $I^*-A$, then we get
	\begin{align*}
	&\frac{d}{dx}\phi^*(W)\\
	&=\frac{d}{dx}\sum_{\sigma\in S_n}\left(\text{sgn}(\sigma)\prod_{i=1}^nc_{i, \sigma(i)}\right)\\
	&=\sum_{\sigma\in S_n}\left(\text{sgn}(\sigma)\frac{d}{dx}\prod_{i=1}^nc_{i, \sigma(i)}\right)\\
	&=\sum_{\sigma\in S_n}\left(\text{sgn}(\sigma)\sum_{j\in\text{fix}(\sigma)}w(j)x^{w(j)-1}\left(\prod_{i\neq j}^nc_{i, \sigma(i)}\right)\right),\\
	\end{align*}
        by the chain rule, where $\text{fix}(\sigma)$ is the set of fixed points of $\sigma$.
	We can switch the order of the summations, but note the permutation must skip the vertex $j$, so it would be a permutation of $\{1, 2, \ldots, j-1, j+1, \ldots, n\}$.
	\begin{align*}
	&=\sum_{j\in V(W)}\sum_{\sigma\in S_{n-1}}\left(\text{sgn}(\sigma)w(j)x^{w(j)-1}\left(\prod_{i\neq j}^nc_{i, \sigma(i)}\right)\right)\\
	&=\sum_{j\in V(W)}w(j)x^{w(j)-1}\phi^*(W\backslash j)
	\end{align*}\qed
\end{enumerate}

\begin{lem}
	Suppose $W_1$, $W_2$ are weighted graphs. Let $W$ be the weighted graph on $n$ vertices obtained by identifying a vertex from $W_1$ and a vertex from $W_2$. Use $v$ to denote the identified vertex, and use $w(v), w_1(v), w_2(v)$ to denote the weight of $v$ in $W$, $W_1$, $W_2$, respectively. Then
	\begin{multline*}
	\phi^*(W)=\phi^*(W_1\backslash v)\phi^*(W_2)+\phi^*(W_1)\phi^*(W_2\backslash v)\\+(x^{w(v)}-x^{w_1(v)}-x^{w_2(v)})\phi^*(W_1\backslash v)\phi^*(W_2\backslash v).
	\end{multline*}
	\label{onesumwr}
\end{lem}

\proof Let $A(W)$ be the adjacency matrix of $W$. We may assume that the $n\times n$ matrix $I^*-A(W)$ is organized so that its first first $\ell$ rows/columns correspond to $W_1$, with the $\ell$-th row/column corresponding $v$, and the remaining rows/columns correspond to vertices in $W_2\backslash v$. Use $a_{ij}$ to denote the $ij$-entry of $I^*-A(W)$, and $M_{i, j}$ to denote the $ij$-cofactor of $I^*-A(W)$. 

We first consider the right hand side of the equation in the lemma statement. Observe $\phi^*(W_1\backslash v)\phi^*(W_2)$ is the product of the determinant of two matrices, namely $I_1^*-A(W_1\backslash v)$ and  $I_2^*-A(W_2)$, where $I_1^*$ and $I_2^*$ are diagonal matrices of appropriate dimensions, whose diagonal entries are $x^{w(i)}$, such that $i$ is the vertex corresponding to the respective row and column. The product $\phi^*(W_1\backslash v)\phi^*(W_2)$ can also be viewed as the determinant of an $n\times n$ matrix $B$ defined as
\[
B:=\begin{bmatrix} I_1^*-A(W_1\backslash v)& \mathbf{0} \\ \mathbf{0} & I_2^*-A(W_2)\end{bmatrix}.
\]
In other words, $B$ is much like the matrix $I^*-A(W)$, with the only differences being that, in $B$, the entries above and to the left of the $\ell\ell$-entry are all zeros, and the $\ell\ell$-entry of $B$ is $x^{w_2(v)}$, while the $\ell\ell$-entry of $I^*-A(W)$ is $x^{w(v)}$. Therefore, if we compute $\phi^*(W_1\backslash v)\phi^*(W_2)$ by applying Laplace's formula to the $\ell$-th column of $B$, the result would be \[\sum_{j>\ell}(-1)^{\ell+j}a_{j\ell}M_{j, \ell}+x^{w_2(v)}M_{\ell, \ell}.\] Note the latter term of the sum is simply $x^{w_2(v)}\phi^*(W_1\backslash v)\phi^*(W_2\backslash v)$.

Similarly, $\phi^*(W_1)\phi^*(W_2\backslash v)$ can be considered as the determinant of a corresponding $n\times n$ matrix as well. So \[\phi^*(W_1)\phi^*(W_2\backslash v)=\sum_{i<\ell}(-1)^{\ell+i}a_{i\ell}M_{i, \ell}+x^{w_1(v)}\phi^*(W_1\backslash v)\phi^*(W_2\backslash v).\]
Likewise
\[
\phi^*(W)=\sum_{k\neq \ell}(-1)^{\ell+k}a_{k\ell}M_{k, \ell}+x^{w(v)}\phi^*(W_1\backslash v)\phi^*(W_2\backslash v).
\]

Therefore,
\begin{multline*}
	\phi^*(W_1\backslash v)\phi^*(W_2)+\phi^*(W_1)\phi^*(W_2\backslash v)\\
	=\phi^*(W)+(x^{w_1(v)}+x^{w_2(v)})\phi^*(W_1\backslash v)\phi^*(W_2\backslash v)-x^{w(v)}\phi^*(W_1\backslash v)\phi^*(W_2\backslash v).
\end{multline*}\qed

With the identities we proved for the weighted characteristic polynomial, we can show that almost all weighted rooted trees have a weighted cospectral mate.

\begin{theorem}
	The proportion of $n$-vertex weighted rooted trees that is weighted cospectral to another weighted rooted tree approaches $1$ as $n$ approaches infinity.
\end{theorem}

\proof Consider the tree $L$ in Figure \ref{schwenkstree}. Let $L_1$ and $L_2$ be the weighted rooted trees obtained by assigning weight $1$ to all vertices in $L$ and assigning $a$ and $b$ as the root, respectively. Observe that $\phi^*(L_1)=\phi^*(L_2)$, and $\phi^*(L_1\backslash a)=\phi^*(L_2\backslash b)$. If $L_1$ is a limb of some weighted rooted tree $W$, then it indicates the weight of $a$ in $W$ is $1$ as well. We could replace $L_1$ in $W$ by $L_2$, and Lemma \ref{onesumwr} implies the resulting weighted rooted tree is weighted cospectral to $W$. By Corollary \ref{haslimbwr}, the proportion of $n$-vertex weighted rooted tress with $L_1$ as a limb approaches $1$ as $n$ approaches infinity. Therefore, the proportion of $n$-vertex weighted rooted trees that is weighted cospectral to another weighted rooted tree approaches $1$ as $n$ approaches infinity.\qed

Moreover, the same proof applies to weighted (unrooted) trees.

\begin{theorem}
	The proportion of $n$-vertex weighted trees that is weighted cospectral to another weighted tree approaches $1$ as $n$ approaches infinity.
\end{theorem}

Thus, we proved that Schwenk's result extends to weighted rooted trees and weighted trees with our definition of the weighted characteristic polynomial.

\section{Closing Remarks}

One main accomplishment in this paper was to extend Schwenk's results to weighted rooted trees. In the case of unweighted trees, the proportion of tree over $n$ vertices with at least one of Schwenk's trees as a limb grows much slower than the proportions of tress with a cospectral mate over $n$ vertices. Looking for constructions for cospectral trees without using Schwenk's limb replacement would help explain this observation.

Since we defined the notion of weighted cospectral for weighted rooted trees, one thing to consider is if there exists any construction for large sets of weighted cospectral weighted rooted trees, other than the limb replacement. It would also be interesting to examine whether the weighted characteristic polynomial would satisfy additional properties exhibited by the characteristic polynomial.

As mentioned at the end of Section 8, the construction we have creates a cut-vertex in the resulting graph, and therefore we could not draw conclusions about general graphs from our construction. To study cospectral vertices in general graphs, it would be useful to look for a construction for $2$-connected graphs with $k$ cospectral vertices for any $k\ge 2$.

One other thing to consider is that when the cospectral vertices resulted from our construction are strongly cospectral. 

Given an $n$-vertex $G$, let $e_i$ denote the standard basis vector in $\mathbb{C}^n$ indexed by the vertex $i\in V(G)$. For each eigenvalue $\theta_r$ of the adjacency matrix $A$ of $G$, there exists an idempotent matrix $E_r$ representing the orthogonal projection onto the eigenspace of the eigenvalue $\theta_r$. Godsil and Smith defined that two vertices $a$ and $b$ in a $n$-vertex graph $G$ are \emph{strongly cospectral} if $E_re_a=\pm E_re_b$ \cite{Godsil1993}. They found that a perfect state transfer from $a$ to $b$ in a quantum walk on $G$ could happen only if $a$ and $b$ are strongly cospectral. Currently, there exists a construction for graphs with pairs of strongly cospectral vertices, but no more than that \cite{GS2017}. The authors of \cite{GS2017} are interested to know if there exists a tree with three pairwise strongly cospectral vertices. We would like to study if a tree resulting from our construction exhibits this property.


\newpage

\bibliographystyle{plain}
\bibliography{refs}


\end{document}